\documentclass[12pt]{article}
\usepackage{epsfig}
\usepackage{color}
\usepackage[english,russian]{babel}
\usepackage[T2A]{fontenc}
\usepackage[utf8]{inputenc}
\usepackage{enumerate,amsmath,amssymb,fullpage}
\usepackage{hyperref}
\usepackage{graphicx, times}

\makeatletter
\def\oversortoftilde#1{\mathop{\vbox{\m@th\ialign{##\crcr\noalign{\kern3\p@}%
      \sortoftildefill\crcr\noalign{\kern3\p@\nointerlineskip}%
      $\hfil\displaystyle{#1}\hfil$\crcr}}}\limits}

\def\sortoftildefill{$\m@th \setbox\z@\hbox{$\braceld$}%
  \braceld\leaders\vrule \@height\ht\z@ \@depth\z@\hfill\braceru$}

\makeatother

\newtheorem{lemma}{Lemma}

\newtheorem{definition}{Definition}
\newtheorem{example}{Example}

\renewcommand{\epsilon}{\varepsilon}

\let\theta\vartheta
\let\phi\varphi

\let\<\langle
\let\>\rangle

\DeclareMathAlphabet{\doba}{U}{msb}{m}{n}

\newcommand{\bq}{\begin{equation}}
\newcommand{\eq}{\end{equation}}

\newcommand{\definedas}{\mathrel{\raise.095ex\hbox{\rm :}\mkern-5.2mu=}}

\makeatletter
\def\thm@headpunct{.}
\makeatother

\begin{document}

\begin{center}
\large \bf SOLVING FREDHOLM INTEGRO-DIFFERENTIAL EQUATIONS USING HYBRID AND BLOCK-PULSE FUNCTIONS 

\end{center}
\vspace{18pt}

\large
\bf \copyright 2025 г. \phantom{20} \  \ \ \ \ \ \ \ \ A. Hosry $^{\it a*}$, R. Nakad $^{\it b**}$, S. Bhalekar$^{\it c***}$ \vspace{12pt}
\normalsize
\rm
\begin{center}\it

$^{\it a}$  Lebanese University, Faculty of Science II, Department of Mathematics, P. O. Box 90656, Fanar, El Metn, Lebanon

$^{\it b}$ Notre Dame University-Louaiz\'e, Faculty of Natural and Applied Sciences, Department of Mathematics and Statistics, P.O. Box 72, Zouk Mikael, Lebanon

$^{\it c}$ School of Mathematics and Statistics, University of Hyderabad,  Hyderabad, 500046 India

e–mail: $^*$ aline.hosry@ul.edu.lb

e–mail: $^{**}$ rnakad@ndu.edu.lb 

e–mail: $^{***}$ sachinbhalekar@uohyd.ac.in, sachin.math@yahoo.co.in

The first two authors are listed by alphabetical order of their family names
\vspace{6pt}

\end{center}

\footnotesize
 \vspace{6pt}
In this paper, hybrid and block-pulse functions are used to approximate the solution of a class of Fredholm integro-differential equations that was first studied by Hemeda \cite{hemada}. By employing suitable approximations, the equation has been converted into a system of algebraic equations that can be solved with classical methods. Finally, the method is explained with illustrative examples and results are compared to the results obtained by Hemeda's method to show the usefulness and efficiency of the block-pulse and hybrid functions approach. 
 \vspace{6pt}

\normalsize


\vspace{12pt}
\centerline{1. INTRODUCTION} 
\vspace{12pt}
Integral equations are central to mathematical models across disciplines like physics, biology, economics and engineering, with many problems reducing to Fredholm integral equations. Researchers have explored their applications in diverse areas, including heat transfer, fluid dynamics, neuroscience, biomechanics, and game theory, among others. 
\\\\
Many authors investigated the general $k^{th}$ order integro-differential equation 
\begin{eqnarray}\label{genequation}
y^{(k)}(t) + l(t) y(t) + \int_a^b g(t, s )y^{(m)}(s) ds = f(t),
\end{eqnarray}
with initial conditions $$y(a) = a_0, \dots, y^{(n-1)}(a) = a_{n-1},$$
where $a_0, \ldots, a_{n-1}$ are real constants, $k, m $ are positive integers with $m< k$, the functions $l, f, g$ are given and $y(t)$ is the solution to be determined.  In \cite{HeHom}, the authors applied the homotopy perturbation method to solve Equation (\ref{genequation}), while in \cite{VIM2, VIM1}, the authors changed the equation to an ordinary integro-differential equation and applied the variational iteration method to solve it so that the Lagrange multipliers can be effectively identified. Using the operational matrix of derivatives of hybrid functions, a numerical method  has been presented in \cite{Hybgen} to solve Equation (\ref{genequation}). In \cite{hemada}, Hemeda used the iterative method introduced in \cite{NIM} to solve the more general  equation 
\begin{eqnarray}\label{genequation2}
y^{(k)}(t) + l(t) y(t) + \int_a^b g(t, s )y^{(n)}(s)y^{(m)}(s) ds = f(t),
\end{eqnarray}
where $n \leq m <k$. 
\\\\
In this paper, we use block-pulse and hybrid functions to approximate solutions $y(t)$ of Equation (\ref{genequation2}). Hybrid functions have been widely utilized to solve differential systems, demon-
strating their usefulness as a valuable mathematical tool. The foundational work in this area was carried out in \cite{B1, B2}, and since then, the hybrid function approach has been refined and extensively applied to approximate differential equations and systems. In this paper, the novelty and the key point in solving Equation (\ref{genequation2}) are to use some useful properties of hybrid functions to derive a new approximation  $Y^{(n)}$ of the derivative  $y^{(n)} (t)$ of order $n$ of the solution $y(t)$ (see Lemma \ref{fund}). Hence,  Equation (\ref{genequation2}) can be converted into  reduced algebraic systems. In Section \ref{Examples}, the method is applied to several examples previously considered by Hemeda, and the results are compared with Hemeda's method to demonstrate the efficiency and effectiveness of the block-pulse and hybrid functions approach.

\vspace{12pt}
\centerline{2. PRELIMINARIES} \label{pre}
\vspace{12pt}

In this section, we define the Legendre polynomials $p_m(t)$, as well as block-pulse and  hybrid functions. We also recall functions' approximation in the Hilbert space $L^2[-1,1]$ (see \cite{Basirat-Malek-Hashemi, Hashemi-Basirat, Hsiao, Malek-Basirat-Hashemi, Sho-Abadi-Golpar}).\\

The Legendre polynomials $p_m(t)$  are polynomials of degree $m$, where $m$ is a positive integer, defined on the interval $[-1,1]$ by
$$p_m(t)= \sum_{k=0}^{M} \frac{(-1)^k (2m-2k)!}{2^m k!(m-k)! (m-2k)!} t^{m-2k}, \;\; m \in \mathbb{N},$$
where $M=   \left \{
\begin{array}{l}
\frac{m}{2}, \qquad  \mbox{if} \; m \; \mbox{is even,} \\ \\
\frac{m-1}{2}, \quad \, \mbox{if} \; m \; \mbox{is odd.}\\
\end{array}
\right.
$\\

Equivalently, the Legendre polynomials are given by the recursive formula 
\[\begin{array}{l}
p_0(t)=1, \quad p_1(t)=t,\\\\
p_{m+1}(t)=\frac{2m+1}{m+1}tp_m(t)-\frac{m}{m+1}p_{m-1}(t), \;\; m=1,2,3,\dots.
\end{array}\]

The set $\{p_m(t), \;m = 0, 1, \dots\}$ is a complete orthogonal system in $L^2[-1, 1]$.

\begin{definition}
For an arbitrary positive integer $q$, let $\{b_k(t)\}_{k=1}^q$ be the finite set of block-pulse functions on the interval $[0, 1)$ defined by
 \begin{eqnarray*}
b_k(t)=\left \{
\begin{array}{l}
1, \quad \mbox{if} \; \frac{k-1}{q} \le t < \frac{k}{q}, \\ \\
0, \quad \mbox{elsewhere.}\\
\end{array}
\right.
\end{eqnarray*}

\end{definition}

The block-pulse functions are disjoint and have the property of orthogonality on $[0, 1)$, since for $i,j=1,2,\dots, q$, we have:
 \begin{eqnarray*}
b_i(t) b_j(t)=\left \{
\begin{array}{l}
0, \qquad \;\mbox{if} \; i \neq j, \\ \\
b_i(t), \quad \mbox{if} \; i = j, \\
\end{array}
\right.
\end{eqnarray*}
and
 \begin{eqnarray*}
\< b_i(t),  b_j(t)\>=\left \{
\begin{array}{l}
0, \quad \mbox{if} \; i \neq j, \\ \\
\frac {1}{q}, \quad \mbox{if} \; i = j,\\
\end{array}
\right.
\end{eqnarray*}
where $\<., .\>$ is the scalar product given by $\<f,g\> = \int_0^1 f(t)g(t)dt$, for any functions $f, g \in L^2[0, 1)$.
\begin{definition}
Let $r$ and $q$ be arbitrary positive integers. The set of hybrid functions $\{b_{km}(t)\},\, k = 1, 2, \dots ,q, \; m = 0, 1, \dots,r- 1$, where $k$ is the order for block-pulse functions, $m$ is the order for Legendre polynomials and $t$ is the normalized time, is defined on the interval $[0, 1)$ as
 \begin{eqnarray*}
b_{km}(t)=\left \{
\begin{array}{l}
p_m(2qt-2k+1), \qquad \mbox{if} \; \frac{k-1}{q} \le t < \frac{k}{q}, \\ \\
0, \qquad \qquad \qquad \qquad \quad \mbox{elsewhere}.\\
\end{array}
\right.
\end{eqnarray*}
\end{definition}

Since Legendre polynomials and block-pulse functions are both complete and orthogonal, then the set of hybrid functions $b_{km}(t)$ is a complete orthogonal system in $L^2[0, 1)$.\\

We are now able to define the $rq \times 1$ vector function $B(t)$ of hybrid functions on $[0, 1)$  by  $$B(t) = \Big (B_1^T (t), \dots, B_q^T (t) \Big)^T,$$
where $B_i (t) = \Big (b_{i0} (t), \dots, b_{i(r-1)}(t)\Big )^T$, for $i=1,2, \dots, q$, and $V^T$ denotes the transpose of a vector $V$. \\

{\bf Function approximation:}  Every function $f(t) \in L^2[0,1)$ can be approximated as
$$f(t) \simeq \sum_{k=1}^{q} \sum_{m=0}^{r-1}  f_{km}b_{km}(t),$$ where
$$f_{km}= \frac{\langle f(t), b_{km}(t)\rangle}{\langle b_{km}(t), b_{km}(t) \rangle}, \, \mbox{for all }  k=1,\dots,q, \; \mbox{for all } m=0,\dots,r-1.$$
Thus,  
\begin{eqnarray}\label{fctF}
f(t) \simeq  F^T B(t) = B^T(t)F,
\end{eqnarray} 
 where $F$ is the $rq \times 1$ column vector having $f_{km}$ as entries.  In a similar way, any function $g(t,s) \in  L^2\big([0, 1) \times [0, 1)\big)$ can be approximated as
\begin{eqnarray}\label{matrixG}
g(t,s) \simeq B^T(t) G B(s),
\end{eqnarray} 
where $G = (g_{ij})$ is the $rq \times rq$ matrix given by
$$g_{ij}= \frac{\langle B_{(i)}(t), \langle g(t,s),B_{(j)}(s)\rangle \rangle}{\langle B_{(i)}(t),B_{(i)}(t)\rangle \langle B_{(j)}(s),B_{(j)}(s)\rangle}, \; i,j=1,2,\dots,rq,$$ and $B_{(i)}(t)$  \big(resp. $B_{(j)}(s)$\big) denotes the $i^{th}$ component \big(resp. the $j^{th}$ component$\big)$ of $B(t)$ \big(resp. $B(s)$\big). \\

{\bf Operational matrix of integration:} The integration of the vector function $B(t)$ may be approximated by $\int_0^t B(s) ds \simeq P B(t)$, where $P$
is an $rq \times rq$ matrix known as the operational matrix of integration and given by

 \[ P =  \left( \begin{array}{ccccc}
E & H & H & ... & H \\
0 & E & H & ... & H \\
0 & 0 & E & ... & ... \\
. & . & . & ...& H \\
0 & 0 & 0 & ... & E
\end{array} \right), \]

where $H$  and $E$ are $r \times r$ matrices defined by

\[  H = \frac 1q  \left( \begin{array}{ccccc}
1 & 0 & 0 & ... & 0 \\
0 & 0 & 0 & ... & 0 \\
0 & 0 & 0 & ... & 0 \\
. & . & . & ...& 0 \\
0 & 0 & 0 & ... & 0
\end{array} \right)\]
 and

\[ E= \frac{1}{2q}\left( \begin{array}{cccccccc}
1 & 1 & 0 & 0 & ... & 0& 0 & 0 \\

-\frac 13 & 0 & \frac 13 & 0 & ... & 0& 0 & 0 \\
0 & -\frac 15 & 0 & \frac 15 & ... & 0& 0 & 0 \\
.. & .. & .. & .. & ... & ..& .. & .. \\
0 & 0 & 0 & 0 & ... & -\frac{1}{2r-3}& 0 & \frac{1}{2r-3} \\
0 & 0 & 0 & 0 & ... & 0& -\frac{1}{2r-1} & 0
\end{array} \right).\] \\\\
{\bf The integration of two hybrid functions:} The integration of the product of two hybrid Legendre block pulse vectors function is given by $L = \int_0^1 B(t) B^T(t) dt,$ where $L$ is the $rq \times rq$ diagonal matrix defined by
 \[ L= \left( \begin{array}{cccccccc}
D & 0 & 0 & 0 & ... & 0& 0 & 0 \\
0 & D  & 0 & 0 & ... & 0 & 0 & 0 \\
.. & .. & .. & .. & ... & .. & .. & .. \\
.. & .. & .. & .. & ... & D & .. & .. \\
0 & 0 & 0 & 0 & ... & 0& D & 0 \\
0 & 0 & 0 & 0 & ... & 0& 0 & D
\end{array} \right)\]

and $D$ is the $r \times r$ matrix given by

 \[ D= \frac 1q \left( \begin{array}{cccccccc}
1 & 0 & 0 & 0 & ... & 0& 0 & 0 \\
0 & \frac 13  & 0 & 0 & ... & 0 & 0 & 0 \\
.. & .. & .. & .. & ... & .. & .. & .. \\
0 & 0 & 0 & 0 & ... & 0& 0 & \frac{1}{2r-1}
\end{array} \right).\] \\

{\bf  The matrix $\widetilde{C}$ associated to a vector $C$:} \label{tilde} For any $rq \times  1$ vector $C$, we define the $rq \times rq$ matrix $\widetilde{C}$ such that
$$B(t)B^T(t) C= \widetilde{C} B(t).$$
$\widetilde{C}$ is called the coefficient matrix. In \cite{Hsiao}, Hsiao computed the matrix $\widetilde{C}$ for $r=2$ and $q=8$, while the authors in \cite{Basirat-Malek-Hashemi} considered the case of $r=4$ and $q=3$.





\vspace{12pt}
\centerline{3. MAIN RESULTS} \label{Main}
\vspace{12pt}

In this section, we approximate solutions $y(t)$ of Equation (\ref{genequation2}). For this, we first start by recalling an approximation for the derivative   $y^{(n)}(t)$, for any $n \geq 1$.
\begin{lemma}\label{fund} \cite{HNB}
Let $y(t)$ be a function and consider its approximation
$y(t)  \simeq Y^{T} B(t) = B^{T}(t) Y$. If $Y^{(n)}$ denotes the approximation of $y^{(n)}(t)$, then for any $n \geq 1$, we have:
$$Y^{(n)} = J^n Y - \sum_{k=1}^n  J^k Y^{(n-k)}_0,$$
where $J =(P^{T})^{-1}$ and $Y_0^{(i)}$ are the approximations of the initial conditions $y_0^{(i)}$, for
$i= 0,\dots, n~-~1$.
\end{lemma}


We are now ready to approximate  Equation (\ref{genequation2}). First, we denote by $V$ the approximation of $l(t)$.  The function $l(t) y(t)$ can be approximated as 
\begin{eqnarray}\label{products}
l(t) y(t) &\simeq& V^T B(t)B^T(t) Y \nonumber \\ &= & \Big(B(t) B^T(t) V\Big)^T Y\nonumber  \\ &\simeq&  \Big( \widetilde V B(t)\Big)^T Y = B^T(t) (\widetilde V)^T Y.
\end{eqnarray}

Using this and  the approximations (\ref{fctF}) and (\ref{matrixG}) of functions of one and two variables, (\ref{genequation2}) can be approximated as
\begin{equation*}
\begin{split}
& B^T(t) Y^{(k)}  +  B^T(t) (\widetilde V)^T Y +  \int_0^ 1 B^T(t) G B(s) B^T(s) Y^{(n)} B^T(s) Y^{(m)} ds = B^T(t) F \\
 & \mbox{which gives }  Y^{(k)}  +(\widetilde V)^T Y +  G \int_0^ 1 B(s) B^T(s) Y^{(n)} B^T(s) Y^{(m)} ds =  F \\
 & \mbox{then we get }  Y^{(k)}  +(\widetilde V)^T Y  +  G \int_0^ 1   \oversortoftilde{Y^{(n)}}B(s)B^T(s) Y^{(m)} ds =  F \\
 & \mbox{ which leads }  Y^{(k)}  +(\widetilde V)^T Y + G  \oversortoftilde{ {Y^{(n)}}}\Big (\int_0^ 1 B(s) B^T(s) ds \Big ) Y^{(m)}= F \\
& \mbox{ and thus, } Y^{(k)}  +(\widetilde V)^T Y  + G \oversortoftilde{ {Y^{(n)}}} L Y^{(m)} =F.
\end{split}
\end{equation*}
Using Lemma \ref{fund}, the last equation becomes
\begin{equation*}
\Big [  J^k Y -\sum_{l=1}^k J^l Y_0^{(k-l)} \Big]  +(\widetilde V)^T Y +  G \Big [  \oversortoftilde{ J^n Y -\sum_{l=1}^n J^l Y_0^{(n-l)}} \Big] L \Big[ {J^m Y -\sum_{l=1}^m J^l Y_0^{(m-l)}}\Big]=F. \label{FredSol}
\end{equation*}

This is a nonlinear system of $rq$ equations in $rq$ variables which can be solved by any iterative method.\\

Wolfram Mathematica is used to solve the algebraic systems arising in this method. We used the commands ``Solve" and ``NSolve" for this purpose. The command ``Solve" provides an exact solution, whenever possible \cite{solve}. On the other hand, the approximate  numerical solution is obtained by using the command ``NSolve" \cite{solve1}.

\vspace{12pt}
\centerline{4. NUMERICAL EXAMPLES}\label{Examples}
\vspace{12pt}

\begin{figure}[h]
	\centering
	\includegraphics[scale=1]{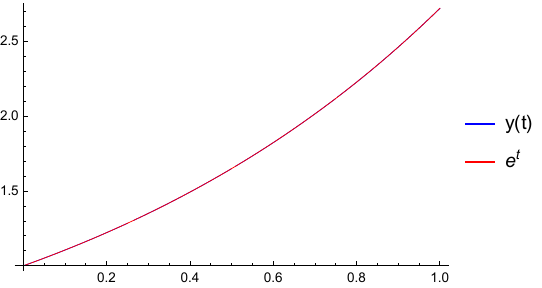}
	\caption{Approximate and exact solutions of System (\ref{ex1}).}
	\label{Fig1}
\end{figure}

\begin{figure}[h]
	\centering
	\includegraphics[scale=1]{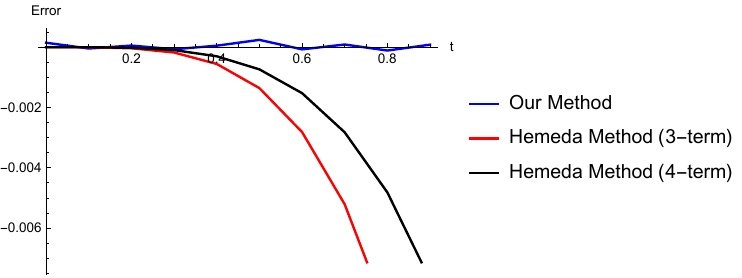}
	\caption{Errors in our method compared with those in Hemeda's three- and four-term approximate solutions.}
	\label{Fig2}
\end{figure}

In this section, we solve some examples and compare our results with those obtained by Hemeda in \cite{hemada}.
\begin{example}
Consider the equation presented in \cite{hemada}
\begin{equation} \label{ex1}
y^{'''} (t) - \int_0^1 st y^{''}(s) y^{''} (s) ds = e^t - \frac {t}{4} (e^2 + 1 )
\end{equation}
with $y(0) = y'(0)= y^{''}(0) = 1$. 
\end{example}
Comparing with Equation (\ref{genequation2}), we get $k=3,$ $m=n=2$, $l(t)=0$, $g(t,s)=-st$ and $f(t)=e^t - \frac {t}{4} (e^2 + 1 )$. We fix $r=3$ and $q=4$. We get
 \begin{eqnarray*}
	B(t)  &=& \left( \chi_{[0, 1/4)},  (-1+8t) \chi_{[0, 1/4)}, (1-24t+96t^2) \chi_{[0, 1/4)}, \chi_{[1/4, 1/2)}, \right.\label{vecB}\\
	&&   (-3+8t) \chi_{[1/4, 1/2)}, (13-72t+96t^2) \chi_{[1/4, 1/2)}, \chi_{[1/2, 3/4)}, (-5+8t)\chi_{[1/2, 3/4)},\nonumber\\
	&& \left. (37-120t+96t^2) \chi_{[1/2, 3/4)}, \chi_{[3/4, 1)}, (-7+8t) \chi_{[3/4, 1)}, (73-168t+96t^2)\chi_{[3/4, 1)}\right), \nonumber
\end{eqnarray*}
(where $\chi_S$ is the characteristic function of a set $S$),
\begin{eqnarray}
	F &=& \left( 0.873944, -0.120293, 0.0059084, 0.672309, -0.0799997, 0.00758654, \right.\nonumber\\
	&&\, \,0.562325, -0.0282622, 0.00974131, 0.570021, 0.0381702, 0.0125081 \left. \right)^T, \nonumber
\end{eqnarray}

\begin{equation*}
	Y_0=Y_1=Y_2 = \left( 1, 0, 0, 1, 0, 0, 1, 0, 0, 1, 0, 0\right)^T, 
\end{equation*}
and
\begin{eqnarray}
	Y &=& \left(1.1361, 0.141865, 0.00591104, 1.45878, 0.182158, 0.00758992, \right.\nonumber\\
	&&\, \, \left. 1.87311, 0.233896, 0.00974566, 2.40513, 0.300328, 0.0125137\right)^T. \nonumber
\end{eqnarray}
This gives the approximate solution of System (\ref{ex1}) as 
\begin{eqnarray}
	y(t) &=& 1.1361 \chi_{[0, 1/4)}+1.45878 \chi_{[1/4, 1/2)} + 1.87312 \chi_{[1/2, 3/4)}+2.40513 \chi_{[3/4, 1)} \nonumber\\
	&&+0.3 (-7+8t) \chi_{[3/4, 1)}+ 0.233896 (-5+8t)\chi_{[1/2, 3/4)} \nonumber\\
	&&  +0.182158 (-3+8t) \chi_{[1/4, 1/2)} + 0.141865(-1+8t) \chi_{[0, 1/4)} \nonumber\\
	&& + 0.0125137(73-168t+96t^2)\chi_{[3/4, 1)}+ 0.00974566(37-120t+96t^2) \chi_{[1/2, 3/4)} \nonumber\\
	&&+ 0.0075899 (13-72t+96t^2) \chi_{[1/4, 1/2)} + 0.00591104 (1-24t+96t^2) \chi_{[0, 1/4)}.\nonumber
\end{eqnarray}
Figure \ref{Fig1} shows the graphs of the approximate solution $y(t)$ and the exact solution $e^t$. The errors at various values of $t$ in our solution (blue color), three-term (red color) and four-term (black color) Hemeda's solution \cite{hemada} are compared in Figure \ref{Fig2}. It can be seen that our method gives better estimate than that in \cite{hemada}.

\begin{figure}[h]
	\centering
	\includegraphics[scale=1]{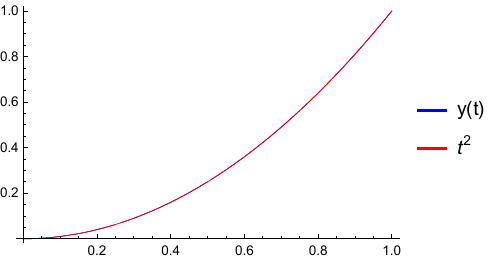}
	\caption{Approximate and exact solutions of System (\ref{ex2}).}
	\label{Fig3}
\end{figure}

\begin{figure}[h]
	\centering
	\includegraphics[scale=1]{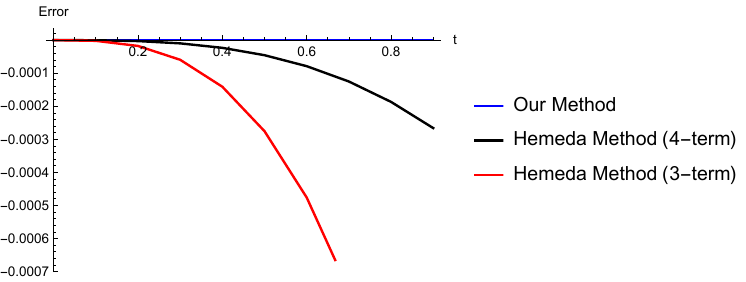}
	\caption{Errors solution of System (\ref{ex2}).}
	\label{Fig4}
\end{figure}

\newpage
\begin{example}
In \cite{hemada}, Hemeda considered the following equation 
\begin{equation}\label{ex2}
y^{''} (t) - \int_0^1 t y(s) y^{'} (s) ds = 2-\frac{t}{2}
\end{equation}
with $y(0) = y'(0)= 0$. 
\end{example}
In fact, we have $k=2, n=0, m=1, l(t)= 0$ and $f(t) =2-\frac{t}{2} $. 
The exact solution of this equation is $y(t) = t^2$. 
For the parameters $r=3$ and $q=4$, $B(t)$ is the same as in the previous example. The other values are 
\begin{eqnarray}
	F &=& \left(1.9375, -0.0625, 0, 1.8125, -0.0625, 8.32667\times10^{-16}, \right.\nonumber\\
	&&\, \,1.6875, -0.0625, -1.33227\times 10^{-14}, 1.5625, -0.0625, -3.55271\times 10^{-14} \left. \right)^T, \nonumber
\end{eqnarray}

\begin{equation*}
	Y_0=Y_1 = \left( 0, 0, 0, 0, 0, 0, 0, 0, 0, 0, 0, 0\right)^T, 
\end{equation*}
and
\begin{eqnarray}
	Y &=& \left(-1.77636\times10^{-15}, 0, 0, -3.55271\times10^{-15}, -2.84217\times10^{-14}, \
	5.68434\times10^{-14}, \right.\nonumber\\
	&&\, \, \left. -1.42109\times10^{-14}, -2.84217\times10^{-14}, -1.13687\times10^{-13}, 1.42109\times10^{-14}, \right.\nonumber\\
	&&\, \, \left.
	2.84217\times10^{-14}, 4.54747\times10^{-13}\right)^T. \nonumber
\end{eqnarray}
This gives the approximate solution of Equation (\ref{ex2}) as 
\begin{eqnarray}
	y(t) &=& 0.0208333 \chi_{[0, 1/4)}+0.145833  \chi_{[1/4, 1/2)} + 0.395833 \chi_{[1/2, 3/4)}+0.770833 \chi_{[3/4, 1)} \nonumber\\
	&&+0.21875 (-7+8t) \chi_{[3/4, 1)}+ 0.15625 (-5+8t)\chi_{[1/2, 3/4)} \nonumber\\
	&&  +0.09375 (-3+8t) \chi_{[1/4, 1/2)} + 0.03125(-1+8t) \chi_{[0, 1/4)} \nonumber\\
	&& + 0.0104167(73-168t+96t^2)\chi_{[3/4, 1)}+ 0.0104167(37-120t+96t^2) \chi_{[1/2, 3/4)} \nonumber\\
	&&+ 0.0104167 (13-72t+96t^2) \chi_{[1/4, 1/2)} + 0.0104167 (1-24t+96t^2) \chi_{[0, 1/4)}.\nonumber
\end{eqnarray}
We sketch this approximate solution along with the exact solution in Figure \ref{Fig3}.
We also compare the errors in our method and Hemeda's method in Figure \ref{Fig4}. As in the previous example, our method works better than the method described in \cite{hemada}. 

\begin{figure}[h]
\centering
\includegraphics[scale=1]{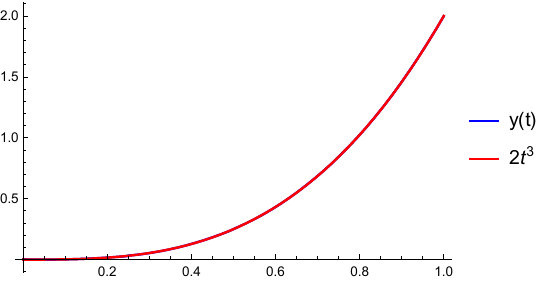}
\caption{Approximate and exact solutions of System (\ref{ex3}).}
\label{Fig5}
\end{figure}

\begin{figure}[h]
\centering
\includegraphics[scale=1]{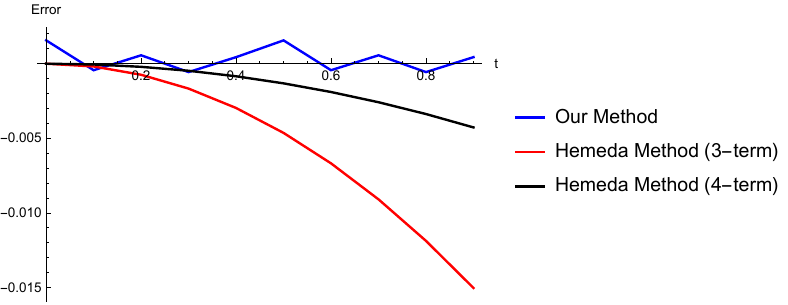}
\caption{Errors in our method compared with those in Hemeda's three- and four-term approximate solutions.}
\label{Fig6}
\end{figure}

\begin{example}
Let us consider the example in \cite{hemada}
\begin{equation}\label{ex3}
	y^{'} (t) - \int_0^1 st y^{2}(s) ds = 6t^2 - \frac {t}{2} 
\end{equation}
and the initial condition $y(0) = 0$. 
\end{example}
In this case, we have $k=1,$ $m=n=0$, $l(t)=0$, $g(t,s)=-st$ and $f(t)= 6t^2 - \frac {t}{2}$. Again, we take $r=3$ and $q=4$. We get

\begin{eqnarray}
F &=& \left(0.0625, 0.125, 0.0625, 0.6875, 0.5, 0.0625, \right.\nonumber\\
&&\, \,2.0625, 0.875, 0.0625, 4.1875, 1.25, 0.0625 \left. \right)^T, \nonumber
\end{eqnarray}
$Y_0$ is the zero vector and 
\begin{eqnarray}
Y &=& \left(6.93889\times10^{-18}, 2.08167\times10^{-17}, 8.32667\times10^{-17}, -2.77556\times10^{-17}, 0, \right.\nonumber\\
&&\, \, \left.
-1.11022\times10^{-15}, 0, -1.11022\times10^{-16}, 3.55271\times10^{-15}, 3.10862\times10^{-15}, \right.\nonumber\\
&&\, \, \left. 4.44089\times10^{-16}, 
-1.77636\times10^{-15}\right)^T. \nonumber
\end{eqnarray}
This gives the approximate solution of (\ref{ex3}) as 
\begin{eqnarray}
y(t) &=& 0.00781248 \chi_{[0, 1/4)}+0.117187 \chi_{[1/4, 1/2)} + 0.507812 \chi_{[1/2, 3/4)}+1.36719 \chi_{[3/4, 1)} \nonumber\\
&&+0.576562 (-7+8t) \chi_{[3/4, 1)}+ 0.295312(-5+8t)\chi_{[1/2, 3/4)} \nonumber\\
&&  +0.107812(-3+8t) \chi_{[1/4, 1/2)} + 0.0140625(-1+8t) \chi_{[0, 1/4)} \nonumber\\
&& + 0.0546875(73-168t+96t^2)\chi_{[3/4, 1)}+ 0.0390625(37-120t+96t^2) \chi_{[1/2, 3/4)} \nonumber\\
&&+ 0.0234375 (13-72t+96t^2) \chi_{[1/4, 1/2)} + 0.00781249  (1-24t+96t^2) \chi_{[0, 1/4)}.\nonumber
\end{eqnarray}
This approximate solution and the exact solution $2t^3$ are shown in Figure \ref{Fig5}. The errors plotted in Figure \ref{Fig6} show the accuracy of our method.

\newpage

\vspace{12pt}
\centerline{5. CONCLUSION}
\vspace{12pt}
In this work, we have introduced a very practical method that involves using hybrid and block-pulse functions to study a Fredholm integro-differential equation that was first solved by Hemeda. The method presented converts the given Fredholm integro-differential equation to a classical algebraic system of equations. The efficiency of the method was illustrated by numerical examples. The figures reported show that the approximate solution is in a good agreement with the exact solution of these examples. We compare the errors in our method and Hemeda's method and show that our method is more powerful. In the future, we can use this method to solve different kinds of integro-differential equations. 

\renewcommand{\refname}{\centerline{REFERENCES}}



\begin{thebibliography}{99}





\bibitem{VIM2} \textit{Abbasbandy S., Shivanian E.} Application of variational iteration method for nth-order integro-differential equations, Zeitschrift fur Naturforschung A, 64 (6-7) (2009), 439-444.

\bibitem{Basirat-Malek-Hashemi} \textit{Basirat B., Maleknejad K., Hashemizadeh E.} Operational matrix approach for the nonlinear Volterra-Fredholm integral equations: Arising in physics and engineering, International Journal of Physical Sciences, 7.2 (2012), 226-233 .




  \bibitem{NIM} \textit{Daftardar-Gejji V., Jafari H.} An iterative method for solving non-linear functional equations, J. Math. Anal. Appl., 316 (2006), 753-763.




\bibitem{HeHom} \textit{Golbabai A.. Javidi M.} Application of He's  homotopy perturbation method for $n^{th}$-order integro-differential equations, Appl. Math. Comput., 190 (2007), 1409-1416.

\bibitem{Hashemi-Basirat} \textit{Hashemizadeh E., Basirat B.} An efficient computational method for the system of linear Volterra integral equations by means of hybrid functions, Mathematical Sciences, 5 (4) (2011), 355-368.

\bibitem{hemada} \textit{Hemeda A. A.} New iterative method: application to $n^{th}$-order integro-differential equations, International Mathematical Forum, Vol. 7 (47) (2012), 2317-2332.



\bibitem{HNB} \textit{Hosry A., Nakad R., Bhalekar S.} A Hybrid function approach to solving a class of Fredholm and Volterra integro-differential equations, Math. Comput. Appl., 2020, 25 (2), 1-16.


\bibitem{Hybgen} \textit{Hou J., Yang C.}, Numerical method in solving Fredholm
integro-differential equations by using hybrid function operational matrix of derivative, Journal of Information and Computational Science, 10:9 (2013), 2757-2764.


\bibitem{Hsiao} \textit{Hsiao C.H.} Hybrid function method for solving Fredholm and Volterra integral equations of the second kind, J. Comput. Appl. Math., 230 (2009), 59-68.


\bibitem{Malek-Basirat-Hashemi} \textit{Maleknejad K., Basirat B., Hashemizadeh E.}, Hybrid Legendre polynomials and block-pulse functions approach for nonlinear Volterra-Fredholm integro-differential equations, Comput. Math. Appl., 61(9) (2011), 282-288.

\bibitem{B1} \textit{Marzban H.R., Razzaghi M.} Optimal control of linear delay systems via
hybrid of block-pulse and Legendre polynomials, J. Franklin Inst., 341(3) (2004), 279-293.

\bibitem{B2}  \textit{Polyanin A.D., Manzhirov A.V.}, Handbook of integral equations (2nd ed.), Chapman and Hall/CRC Press, Boca Raton- London (2008).







\bibitem{VIM1} \textit{Shang X., Han D.}, Application of the variational iteration method for solving $n^{th}$-order integro-differential equations, Journal of Computational and Applied Mathematics, 234 (5) (2010), 1442-1447.


\bibitem{Sho-Abadi-Golpar} \textit{Shojaeizadeh T., Abadi Z., Golpar Raboky E.} Hybrid functions approach for
solving Fredholm and Volterra integral equations, J. Prime Res. Math., 5 (2009), 124-132.

\bibitem {solve} Wolfram Research (1988), Solve, Wolfram Language function, https://reference.wolfram.com/language/ref/Solve.html (updated 2024).
\bibitem {solve1} Wolfram Research (1991), NSolve, Wolfram Language function, https://reference.wolfram.com/language/ref/NSolve.html (updated 2024).

\end{thebibliography}
\end{document}